\documentclass[12pt]{article}
\usepackage{mathrsfs}
\usepackage{amsmath,amsfonts,amssymb,rotating,amsthm}
\usepackage{hyperref}
\usepackage{array}
\usepackage{color}

\def\qed{\nopagebreak\hfill{\rule{4pt}{7pt}}}
\def\proof{\noindent {\it{Proof.} \hskip 2pt}}

\newtheorem{theo}{Theorem}[section]

\newtheorem{lemm}[theo]{Lemma}

\newtheorem{conj}[theo]{Conjecture}

\theoremstyle{remark}

\newdimen\Squaresize \Squaresize=11pt
\newdimen\Thickness \Thickness=0.7pt
\def\Square#1{\hbox{\vrule width \Thickness
   \vbox to \Squaresize{\hrule height \Thickness\vss
    \hbox to \Squaresize{\hss#1\hss}
   \vss\hrule height\Thickness}
\unskip\vrule width \Thickness} \kern-\Thickness}

\def\Vsquare#1{\vbox{\Square{$#1$}}\kern-\Thickness}

\def\moins{\raise 1pt\hbox{{$\scriptstyle -$}}}
\begin{document}
\begin{center}
%\bf{NOTE}
\end{center}

\begin{center}
{\large \bf Two-log-convexity of the Catalan-Larcombe-French sequence}
\end{center}
\begin{center}
Brian Yi Sun$^1$, Baoyindureng Wu$^2$\\[6pt]

College of Mathematics and System Science,\\
XinJiang University, Urimqi, 830046, P. R. China

Email: $^{1}${\tt brianys1984@126.com},
       $^{2}${\tt wubaoyin@hotmail.com}
\end{center}

\vspace{0.3cm} \noindent{\bf Abstract.}  The Catalan-Larcombe-French sequence $\{P_n\}_{n\geq 0}$ arises in a series expansion of the complete elliptic integral of the first kind. It has been
proved that the sequence is log-balanced. In the paper, by exploring
a criterion due to Chen and Xia for testing 2-log-convexity of a sequence
satisfying three-term recurrence relation, we prove that the new sequence
$\{P^2_n-P_{n-1}P_{n+1}\}_{n\geq 1}$ are strictly log-convex and hence the Catalan-Larcombe-French sequence is strictly 2-log-convex.



\section{Introduction}

This paper is concerned with the log-behavior of the Catalan-Larcombe-French sequence.
To begin with, let us recall that a sequence $\{z_n\}_{n\geq 0}$ is said to be log-concave if
\begin{equation}\label{def-log-concave}
z_n^2\geq z_{n+1}z_{n-1},\,\,\,for\,\,\, n\geq 1,
\end{equation}
and it is log-convex if
\begin{equation}\label{def-log-convex}
z_n^2\leq z_{n+1}z_{n-1},\,\,\,for\,\,\, n\geq 1.
\end{equation}
Meanwhile, the sequence $\{z_n\}_{n\geq 0}$ is called strictly log-concave{(resp. log-convex)} if the inequality in \eqref{def-log-concave}{(resp. \eqref{def-log-convex})} is strict for all
$n\geq 1.$
We call $\{z_n\}_{n\geq 0}$  log-balanced if the sequence itself is log-convex while $\{\frac{z_n}{n!}\}_{n\geq 0}$ is log-concave.

Given a sequence $A=\{z_n\}_{n\geq 0}$, define the operator $\mathcal{L}$ by
$$\mathcal{L}(A)=\{s_n\}_{n\geq 0},$$
where $s_n=z_{n-1}z_{n+1}-z_n^2$ for $n\geq 1.$ We say that $\{z_n\}_{n\geq 0}$
is $k$-log-convex({\it{resp. $k$-log-concave}}) if $\mathcal{L}^j(A)$ is log-convex({\it{resp. log-concave}}) for all $j=0,1,\cdots,k-1$, and that $A=\{z_n\}_{n\geq 0}$ is $\infty$-log-convex({\it{resp. $\infty$-log-concave}}) if $\mathcal{L}^k(A)$ is
log-convex({\it{resp. log-concave}}) for any $k\geq 0.$ Similarly, we can define strict $k$-log-concavity or strict $k$-log-convexity of a sequence.

It is worthy to mention that besides that they are fertile sources
of inequalities, log-convexity and log-concavity have many applications in some different mathematical disciplines, such as geometry, probability theory, combinatorics and so on. See the surveys due to Brenti \cite{brenti} and Stanley \cite{stanley} for more details. Additionally, it is clear that the log-balancedness implies the log-convexity and a sequence $\{z_n\}_{n\geq 0}$ is log-convex(resp.\,\, log-concave) if and only if its quotient sequence $\{\frac{z_n}{z_{n-1}}\}_{n\geq 1}$ is nondecreasing(resp.\,\,nonincreasing).
It is also known that the quotient sequence of a log-balanced sequence does not grow too fast. Therefore, log-behavior are important properties of combinatorial sequences and they are instrumental in obtaining the growth rate of a sequence. Hence the log-behaviors of a sequence deserves to be investigated.

In this paper, we investigate the 2-log-behavior of the Catalan-Larcombe-French sequence, denoted by $\{P_n\}_{n\geq 0}$,  which arises in connection with series expansions of the complete elliptic integrals of the first kind \cite{jv2010,zhao2014}. To be precise, for $0<|c|<1$,
\begin{equation*}
\int_0^{\pi/2}\frac{1}{\sqrt{1-c^2\sin^2{\theta}}}d\theta=\frac{\pi}{2}
\sum_{n=0}^\infty\left(\frac{1-\sqrt{1-c^2}}{16}\right)^nP_n.
\end{equation*}

Furthermore, the numbers $P_n$ can be  written as the following sum:
\begin{align*}%\label{CLF-sum}
P_n=2^n\sum_{k=0}^{n} (-4)^k\binom{n-k}{k}\binom{2n-2k}{n-k}^2,
\end{align*}
see \cite[A05317]{Sloane}.
Besides,  the number $P_n$  satisfies three-term recurrence relations \cite{zhao2014} as follows:
\begin{align}\label{C-L-F_recurrence}
(n+1)^2 P_{n+1}=8(3n^2+3n+1)P_n-128n^2P_{n-1}, \,\,for \quad n\geq 1,
\end{align}
with the initial values $P_0=1$ and $P_1=8$.

Recently, Zhao \cite{zhao2014} studied the log-behavior of the Catalan-Larcombe-French sequence and proved that the sequence $\{P_n\}_{n\geq 0}$ is log-balanced. What's more, the Catalan-Larcombe-French sequence has many interesting properties and the reader can refer  \cite{jv2010,lf2000,zhao2014}. In the sequel, we study the 2-log-behavior of the sequences and obtain the following result.
\begin{theo}\label{2-log-convex-CLF}
The Catalan-Larcombe-French sequence $\{P_n\}_{n\geq 0}$ is strictly 2-log-convex, that is,
\begin{equation}
\mathcal{P}_n^2<\mathcal{P}_{n-1}\mathcal{P}_{n+1},
\end{equation}
where $\mathcal{P}_n=P_n^2-P_{n-1}P_{n+1}.$
\end{theo}
 We will give our proof of Theorem \ref{2-log-convex-CLF} in the third section by utilizing a testing criterion, which is proposed by Chen and Xia \cite{ChenXia}.

 To make this paper self-contained, let us recall their criterion.

\begin{theo}[Chen and Xia\cite{ChenXia}]\label{chenxia}
Suppose $\{z_n\}_{n\geq 0}$ is a positive log-convex sequence that satisfies the
following three-term recurrence relation
\begin{equation}\label{three-term recurrence}
z_n=a(n)z_{n-1}+b(n)z_{n-2}, \,\,for\,\,\, n\geq 2.
\end{equation}
Let
\begin{align*}
c_0(n)=&-b^2(n+1)[a^2(n+2)+b(n+1)-a(n+2)a(n+3)-b(n+3)];\\
c_1(n)=&b(n+1)[2a(n+2)b(n+1)+2a(n+3)a(n+2)a(n+1)\\
&+a(n+3)b(n+2)+2a(n+1)b(n+3)-2a^2(n+2)a(n+1)\\
&-2a(n+2)b(n+2)-3a(n+1)b(n+1)];\\
c_2(n)=&4a(n+1)a(n+2)b(n+1)+2b(n+1)b(n+2)+a^2(n+1)a(n+2)a(n+3)\\
&+a(n+1)a(n+3)b(n+2)+a^2(n+1)b(n+3)-3a^2(n+1)b(n+1)\\
&-a(n+3)a(n+2)b(n+1)-a^2(n+2)a^2(n+1)-b(n+3)b(n+1)\\
&-2a(n+2)a(n+1)b(n+2)-b^2(n+2);\\
c_3(n)=&2a^2(n+1)a(n+2)+2a(n+1)b(n+2)-a(n+1)b(n+3)-a^3(n+1)\\
&-a(n+1)a(n+2)a(n+3)-a(n+3)b(n+2);
\end{align*}
and
$$\Delta(n)=4c_2^2(n)-12c_1(n)c_3(n).$$
Assume that $c_3(n)<0$ and $\Delta(n)\geq 0$ for all $n\geq N$, where
$N$ is a positive integer. If there exist $f_n$ and $g_n$ such that, for all
$n\geq N$,
\begin{itemize}
\item[(I)]$f_n\leq \frac{z_n}{z_{n-1}}\leq g_n;$
\item[(II)]$f_n\geq \frac{-2c_2(n)-\sqrt{\Delta(n)}}{6c_3(n)};$
\item[(III)]$c_3(n)g_n^3+c_2(n)g_n^2+c_1(n)g_n+c_0(n)\geq 0,$
\end{itemize}
then we see that
$\{z_n\}_{n\geq N}$ is  $2$-log-convex,
that is, for $n\geq N$,
\begin{equation*}
(z_{n-1}z_{n+1}-z_n^2)(z_{n+1}z_{n+3}-z_{n+2}^2)>(z_{n}z_{n+2}-z_{n+1}^2)^2.
\end{equation*}
 \end{theo}

With respect to the theory in this field, it should be mentioned that the log-behavior of a sequence which satisfies a three-term recurrence has been extensively studied; see Liu and Wang\cite{liuwang}, Chen et al. \cite{CGW2014,CGW}, Liggett \cite{Liggett}, Do\v{s}li\'{c}\cite{Doslic}, etc.

\section{Bounds for $\frac{P_n}{P_{n-1}}$}\label{bunds-sec}

Before proving  Theorem \ref{2-log-convex-CLF}, we need the following two lemmas.
\begin{lemm}\label{lower-bound-ratio-CLF}
Let $$f_n=\frac{232n}{15(n+2)},$$ and $P_n$ be the sequence defined by the recurrence relation \eqref{C-L-F_recurrence}. Then we have, for all $n\geq 1,$
\begin{equation}\label{lower-bound-inequality}
\frac{P_n}{P_{n-1}}> f_n.
\end{equation}
\end{lemm}
\proof We proceed the proof by induction. First note that, for $n=1$ and $n=2$, we have
$\frac{P_1}{P_0}=8>\frac{232}{45}$ and $\frac{P_2}{P_1}=10>\frac{464}{60}$.
Assume that the inequality \eqref{lower-bound-inequality} is valid for $n\leq k.$
We will show that
$$\frac{P_{k+1}}{P_{k}}> f_{k+1}.$$
By the recurrence \eqref{C-L-F_recurrence}, we have
\begin{align*}
\frac{P_{k+1}}{P_k}&=\frac{8(3k^2+3k+1)}{(k+1)^2}-
\frac{128k^2}{(k+1)^2}\frac{P_{k-1}}{P_k}\\
&>\frac{8(3k^2+3k+1)}{(k+1)^2}-\frac{128k^2}{(k+1)^2}\frac{1}{f_k}\\
&=\frac{8(57k^2+27k+29)}{29(k+1)^2}\\
&>f_{k+1},
\end{align*}
in which the last inequality follows by
\begin{align*}
\frac{8(57k^2+27k+29)}{29(k+1)^2}-f_{k+1}=\frac{8(14k^3+447k^2-873k+464)}{435(k+1)^2(k+3)}>0,
\end{align*}
for all $k\geq 1.$
This completes the proof.
\qed

\begin{lemm}\label{upper-bound-ratio-CLF}
Let $$g_n=16-\frac{16}{n}-\frac{16}{n^3},$$ and $P_n$ be the sequence defined by the recurrence relation \eqref{C-L-F_recurrence}. Then we have, for all $n\geq 6,$
\begin{equation}\label{upper-bound-inequality}
\frac{P_n}{P_{n-1}}\leq g_n.
\end{equation}
\end{lemm}
\proof
 First note that, for $n=6$, we have
$\frac{P_6}{P_5}=\frac{3562}{269}<g_6=\frac{358}{27}$.
Assume that for $k\geq 6$, the inequality \eqref{upper-bound-inequality} is valid for $n\leq k.$
We will show that
$$\frac{P_{k+1}}{P_{k}}< g_{k+1}.$$
By the recurrence \eqref{C-L-F_recurrence}, we have
\begin{equation}\label{inequality_1}
\begin{split}
\frac{P_{k+1}}{P_k}&=\frac{8(3k^2+3k+1)}{(k+1)^2}-
\frac{128k^2}{(k+1)^2}\frac{P_{k-1}}{P_k}\\
&<\frac{8(3k^2+3k+1)}{(k+1)^2}-\frac{128k^2}{(k+1)^2}\frac{1}{g_k}\\
&=\frac{8(2k^5-2k^3-4k^2-3k-1)}{(k+1)^2(k^3-k^2-1)}.
\end{split}
\end{equation}
Consider
\begin{align}
\frac{8(2k^5-2k^3-4k^2-3k-1)}{(k+1)^2(k^3-k^2-1)}-g_{k+1}=
-\frac{8(5k^2+2k+3)}{(k+1)^3(k^3-k^2-1)}<0,\label{inequality_2}
\end{align}
for all $k\geq 2.$
So we see that for all $n\geq 6$, the inequality \eqref{upper-bound-inequality} holds by induction.

\qed

 With the above lemmas in hand, we are now in a position to prove our main result in the next section.

\section{Proof of  Theorem \ref{2-log-convex-CLF}}\label{proof-sec}
In this section, by using the criterion of Theorem \ref{chenxia}, we can show that the Catalan-Larcombe-French sequence is strictly 2-log-convex.

To begin with, the following lemma, which is obtained by Zhao \cite{zhao2014}, is indispensable for us.
\begin{lemm}[Zhao\cite{zhao2014}]\label{zhao}
The Catalan-Larcombe-French sequence is log-balanced.
\end{lemm}
By the definition of log-balanced sequence, we know that $\{P_n\}_{n\geq 0}$ is
log-convex.

{\it Proof of Theorem \ref{2-log-convex-CLF}.}~
By {\it Lemma \ref{zhao}}, it suffices for us to show that
\begin{equation*}
(P_{n-1}P_{n+1}-P_n^2)(P_{n+1}P_{n+3}-P_{n+2}^2)-(P_{n}P_{n+2}-P_{n+1}^2)^2>0.
\end{equation*}
According to the recurrence relation \eqref{C-L-F_recurrence}, we see that
\begin{align*}
a(n)=&\frac{8(3n^2-3n+1)}{n^2};\\
b(n)=&-\frac{128(n-1)^2}{n^2}.
\end{align*}
By taking $a(n),b(n)$ in $c_0,\cdots,c_3$, we can obtain
\begin{align*}
c_3(n)=&-\frac{512}{(n+1)^6(n+2)^2(n+3)^2}\left(3n^8+5n^7-27n^6-32n^5+112n^4\right.\\
&\left.+234n^3+177n^2+63n+9\right)\\
&<0,
\end{align*}
for all $n\geq 1.$
Besides, we have to verify that for some positive integer $N$, the conditions $(II)$ and $(III)$ in {\it Theorem \ref{chenxia}} hold for all
$n\geq N$. That is,
\begin{align}
f_n\geq \frac{-2c_2(n)-\sqrt{\Delta(n)}}{6c_3(n)};\label{condition-I}\\
c_3(n)g_n^3+c_2(n)g_n^2+c_1(n)g_n+c_0(n)\geq 0.\label{condition-II}
\end{align}
Let $$\delta(n)=-6c_3(n)f_n-2c_2(n)$$ and $$f(g_n)=c_3(n)g_n^3+c_2(n)g_n^2+c_1(n)g_n+c_0(n).$$
To show \eqref{condition-I}, it is equivalent to show that, for some positive integers $N$,
$\delta(n)\geq 0$ and $\delta^2(n)\geq \Delta(n).$
By calculating, we easily find that, for all $n\geq 1,$
\begin{align*}
\delta(n)=&\frac{8192}{5(n+1)^6(n+2)^4(n+3)^2}\left(32n^{10}+129n^9+472n^8+3556n^7+12157n^6\right.\\
&\left.+17632n^5+10550n^4+1293n^3-1500n^2-798n-135\right)\\
&\geq 0,
\end{align*}
and for all $n\geq 3,$
\begin{align*}
\delta^2(n)-\Delta(n)=&\frac{67108864n}{25(n+3)^4(n+2)^7(n+1)^{12}}
\left(699n^{18}+2158n^{17}+6983n^{16}\right.\\
&+97994n^{15}+155517n^{14}-1256916n^{13}-3302168n^{12}+\\
&5191280n^{11}+25505142 n^{10}+14486584 n^9-63005002 n^8\\
&-153766236 n^7-178037517 n^6-131841558n^5-68012397n^4\\
&\left.-24910146n^3-6269211n^2-975888n-70470\right)\\
&\geq 0.
\end{align*}
Thus, take $N=3$ and for all $n\geq N$, we have
$\delta(n)\geq 0,\delta^2(n)\geq \Delta(n),$ which follows from
the inequality \eqref{condition-I}. We show the inequality \eqref{condition-II} for some positive integer $M$. Note that, by {\it Lemma \ref{upper-bound-ratio-CLF} } and some calculations, we have
\begin{align*}
f(g_n)&=c_3(n)g_n^3+c_2(n)g_n^2+c_1(n)g_n+c_0(n)\\
=&\frac{1048576}{n^9(n+1)^6(n+2)^4(n+3)^2}
\left(54n^{15}+378 n^{14}+916 n^{13}+644 n^{12}-1529 n^{11}\right.\\
&-5340 n^{10}-8383 n^9-7416n^8-2284n^7+4156n^6+7969n^5+7688n^4\\
&\left.+4953n^3+2154n^2+576n+72\right).
\end{align*}
Take $M=6$, it is not difficult to verify that for all $n\geq M$,
$$f(g_n)>0.$$
Let $N_0=max\{N,M\}=6$, then for all $n\geq 6$, all of the above inequalities hold.
By {\it Lemma \ref{zhao} and Theorem \ref{chenxia}}, the Catalan-Larcombe-French sequence $\{P_n\}_{n\geq 6}$ is strictly 2-log-convex for all $n\geq 6.$
What is more, one can easily test that these numbers $\{P_n\}_{0\leq n\leq 8}$ also satisfy the property of 2-log-convexity by simple calculations.
Therefore, the whole sequence $\{P_n\}_{n\geq 0}$ is strictly 2-log-convex.
This completes the proof.
\qed

It deserves to be mentioned that by considerable calculations and plenty of verifications, the following conjectures should be true.
\begin{conj}
The Catalan-Larcombe-French sequence is $\infty$-log-convex.
\end{conj}
\begin{conj}
The quotient sequence $\{\frac{P_n}{P_{n-1}}\}_{n\geq 1}$ of the Catalan-Larcombe-French sequence is log-concave, equivalently,
for all $n\geq 2,$
$$P_{n-2}P_n^3\geq P_{n+1}P_{n-1}^3.$$
\end{conj}

\noindent{\bf Acknowledgments.} This work was supported Research supported by NSFC (No. 11161046) and by the Xingjiang Talent Youth Project (No. 2013721012).

\end{document}